\documentstyle[amscd,amssymb,verbatim,12pt]{amsart}

\theoremstyle{plain}
\newtheorem{Thm}[subsection]{Theorem}
\newtheorem{Cor}[subsection]{Corollary}
\newtheorem{Lem}[subsection]{Lemma}
\newtheorem{Prop}[subsection]{Proposition}

\theoremstyle{definition}
\newtheorem{Def}[subsection]{Definition}

\theoremstyle{remark}

\newtheorem{Rem}[subsection]{Remark}

\numberwithin{equation}{section}
\renewcommand{\rm}{\normalshape}

\newif\ifShowLabels
\ShowLabelstrue
\newdimen\theight
\def\TeXref#1{%
	\leavevmode\vadjust{\setbox0=\hbox{{\tt
		\quad\quad  {\small \rm #1}}}%
	\theight=\ht0
	\advance\theight by \lineskip
	\kern -\theight \vbox to
	\theight{\rightline{\rlap{\box0}}%
	\vss}%
	}}%

\renewcommand{\sec}[2]{\section{#2}\label{S:#1}%
	\ifShowLabels \TeXref{{S:#1}} \fi}
\newcommand{\ssec}[2]{\subsection{#2}\label{SS:#1}%
	\ifShowLabels \TeXref{{SS:#1}} \fi}

\newcommand{\refs}[1]{Section ~\ref{S:#1}}
\newcommand{\refss}[1]{Section ~\ref{SS:#1}}
\newcommand{\reft}[1]{Theorem ~\ref{T:#1}}
\newcommand{\refl}[1]{Lemma ~\ref{L:#1}}
\newcommand{\refp}[1]{Proposition ~\ref{P:#1}}

\newcommand{\refd}[1]{Definition ~\ref{D:#1}}
\newcommand{\refr}[1]{Remark ~\ref{R:#1}}
\newcommand{\refe}[1]{\eqref{E:#1}}

\newenvironment{thm}[1]%
	{ \begin{Thm} \label{T:#1}  \ifShowLabels \TeXref{T:#1} \fi }%
	{ \end{Thm} }

\newcommand{\th}[1]{\begin{thm}{#1} }
\renewcommand{\eth}{\end{thm} }

\newenvironment{lemma}[1]%
	{ \begin{Lem} \label{L:#1}  \ifShowLabels \TeXref{L:#1} \fi }%
	{ \end{Lem} }
\newcommand{\lem}[1]{\begin{lemma}{#1}}
\newcommand{\elem}{\end{lemma}}

\newenvironment{propos}[1]%
	{ \begin{Prop} \label{P:#1}  \ifShowLabels \TeXref{P:#1} \fi }%
	{ \end{Prop} }
\newcommand{\prop}[1]{\begin{propos}{#1}}
\newcommand{\eprop}{\end{propos}}

\newenvironment{corol}[1]%
	{ \begin{Cor} \label{C:#1}  \ifShowLabels \TeXref{C:#1} \fi }%
	{ \end{Cor} }
\newcommand{\cor}[1]{\begin{corol}{#1}}
\newcommand{\ecor}{\end{corol}}

\newenvironment{defeni}[1]%
	{ \begin{Def} \label{D:#1}  \ifShowLabels \TeXref{D:#1} \fi }%
	{ \end{Def} }
\newcommand{\defe}[1]{\begin{defeni}{#1}}
\newcommand{\edefe}{\end{defeni}}

\newenvironment{remark}[1]%
	{ \begin{Rem} \label{R:#1}  \ifShowLabels \TeXref{R:#1} \fi }%
	{ \end{Rem} }
\newcommand{\rem}[1]{\begin{remark}{#1}}
\newcommand{\erem}{\end{remark}}

\newcommand{\eq}[1]%
	{ \ifShowLabels \TeXref{E:#1} \fi
	   \begin{equation} \label{E:#1} }
\newcommand{\eeq}{ \end{equation} }

\newcommand{\prf}{ \begin{pf} }
\newcommand{\epr}{ \end{pf} }

\newcommand\alp{\alpha}		
\newcommand\bet{\beta}
		
\newcommand\del{\delta}		\newcommand\Del{\Delta}
\newcommand\eps{\varepsilon}
\newcommand\zet{\zeta}

\newcommand\lam{\lambda}		
\newcommand\sig{\sigma}		

\newcommand\ome{\omega}		\newcommand\Ome{\Omega}

\newcommand\calB{{\cal{B}}}

\newcommand\calE{{\cal{E}}}
\newcommand\calF{{\cal{F}}}

\newcommand\calM{{\cal{M}}}

\newcommand\calO{{\cal{O}}}

\newcommand\calV{{\cal{V}}}

\newcommand\calZ{{\cal{Z}}}

\newcommand\RR{\Bbb{R}}

\newcommand\ZZ{\Bbb{Z}}

\newcommand\CC{\Bbb{C}}

\newcommand\NN{\Bbb{N}}

\newcommand\nek{,\ldots,}
\newcommand\sdp{\times \hskip -0.3em {\raise 0.3ex
\hbox{$\scriptscriptstyle |$}}} 

\newcommand\End{\operatorname{End\,}}

\newcommand\Ker{\operatorname{Ker}}

\newcommand\rk{\operatorname{rk}}

\newcommand\RE{\operatorname{Re}}

\newcommand\SPAN{\operatorname{span}}

\newcommand\Tr{\operatorname{Tr}}

\newcommand\od{{\overline{d}}}

\newcommand\om{{\overline{m}}}

\newcommand\olam{{\overline{\lambda}}}

\newcommand\tilV{{\widetilde{V}}}

\newcommand\tilDel{{\widetilde{\Delta}}}

\newcommand\tilpat{{\widetilde{\partial}}}

\newcommand\ten{\otimes}

\ShowLabelsfalse

\theoremstyle{plain}

\renewcommand{\b}{\bullet}
\newcommand{\h}[1]{\text{\( H^{#1}(M,F)\)}}
\newcommand{\hb}{\h{\b}}
\newcommand{\dh}[1]{\text{\( \det H^{#1}(M,F)\)}}
\newcommand{\dhb}{\text{\( \det H^{\b}(M,F)\)}}

\newcommand{\ts}{\text{\( C^{\b}(W^u,F)\)}}
\newcommand{\dts}{\text{\( \det C^{\b}(W^u,F)\)}}

\renewcommand{\om}{\text{\( \Ome^{\b}(M,F)\)}}

\newcommand{\omot}{\text{\( \Ome^{\b,[0,1]}_t(M,F)\)}}

\newcommand{\fdm}{\text{\( \|\cdot\|_{\det H^{\b}(V)}\)}}
\newcommand{\fdmb}{\text{\( |\cdot|_{\det H^{\b}(V)}\)}}

\newcommand{\rsm}{\text{\( \|\cdot\|^{RS}_{\det \h{\b}}\)}}
\newcommand{\rsmt}{\text{\( \|\cdot\|^{RS}_{\det \h{\b},t}\)}}

\newcommand{\mm}{\text{\( \|\cdot\|^\calM_{\det \h{\b}}\)}}
\newcommand{\mmt}{\text{\( \|\cdot\|^\calM_{\det \h{\b},t}\)}}

\newcommand{\rsmb}{\text{\( |\cdot|^{RS}_{\det \h{\b}}\)}}
\newcommand{\rsmbt}{\text{\( |\cdot|^{RS}_{\det \h{\b},t}\)}}

\newcommand{\mmb}{\text{\( |\cdot|^\calM_{\det \h{\b}}\)}}
\newcommand{\mmbt}{\text{\( |\cdot|^\calM_{\det \h{\b},t}\)}}

\newcommand{\mt}{\text{\( \rho^\calM\)}}
\newcommand{\mtt}{\text{\( \rho^\calM(t)\)}}

\newcommand{\rst}{\text{\( \rho^{RS}\)}}
\newcommand{\rstt}{\text{\( \rho^{RS}(t)\)}}
\newcommand{\sst}{\text{\( \rho^{ss}\)}}

\newcommand{\gtm}{\text{\( g^{TM}\)}}
\newcommand{\gft}{\text{\( g^F_t\)}}
\newcommand{\gf}{\text{\( g^F\)}}

\newcommand{\nf}{\text{\( \nabla f\)}}
\newcommand{\ntm}{\text{\( \nabla^{TM}\)}}
\newcommand{\etm}{\text{\( e(TM,\ntm)\)}}
\newcommand{\mqf}{\text{\( \psi(TM,\ntm)\)}}
\newcommand{\thet}{\text{\( \theta(F,\gf)\)}}

\newcommand{\p}{\text{\( \partial\)}}
\newcommand{\tp}{\text{\(\widetilde \partial\)}}

\begin{document}

\title[Witten deformation of the analytic torsion]
	{Witten deformation of the analytic torsion and \\
	the spectral sequence of a filtration}
\author{Maxim Braverman}
\address{School of Mathematical Sciences\\
Tel-Aviv University\\
Ramat-Aviv 69978, Israel}
\email{maxim@@math.tau.ac.il}

\thanks{The research was supported by grant No. 449/94-1 from the
Israel Academy of Sciences and Humanities.}


\maketitle

\begin{abstract}
 Let $F$ be a flat vector bundle over a compact Riemannian manifold $M$
and let $f:M\to \RR$ be a Morse function.
Let $g^F$ be a smooth Euclidean
metric on $F$, let $g^F_t=e^{-2tf}g^F$  and let $\rho^{RS}(t)$
be the Ray-Singer analytic torsion of $F$ associated to the metric
$g^F_t$. Assuming that $\nabla f$ satisfies the Morse-Smale
transversality conditions,
we provide an asymptotic expansion for $\log \rho^{RS}(t)$
for $t\to +\infty$
of the form $a_0+a_1t+b\log\left(\frac t\pi\right)+o(1)$,
where the coefficient $b$ is a half-integer depending only on the
Betti numbers of $F$.
In the case where all the critical values of $f$ are rational, we
calculate the coefficients $a_0$ and $a_1$  explicitly
in terms of the spectral sequence of a filtration associated to the
Morse function.
These results are obtained as an applications of a theorem by Bismut
and Zhang.
\end{abstract}

\sec{introd} {Introduction}

\ssec{rst}{} The {\em analytic torsion},
\rst, introduced by Ray and Singer
\cite{rs}, is a numerical invariant associated to a flat Euclidean
vector bundle $F$ over a compact Riemannian manifold $M$.
It depends smoothly on the Riemannian metric  \gtm\ on $TM$
and on the Euclidean metric \gf\ on $F$.

Let $f:M\to\RR$ be a Morse function.  By Witten deformation of the
Euclidean bundle $F$ we shall understand the family of metrics
\eq{gft}
   	\gft=e^{-2tf}\gf, \qquad    t>0
\end{equation}
on $F$. Let \rstt\ be the Ray-Singer torsion of $F$ associated to
the metrics \gtm\ and \gft.

Burghelea, Friedlander and Kappeler
(\cite{bfk3}) have shown  (with some additional conditions on $f$, cf.
\refss{assum}) that  the function $\log\rstt$ has an
asymptotic expansion for $t\to +\infty$ of the form
\eq{bfk}
	\log\rst(t)=\sum_{j=0}^{n+1}a_jt^j+b\log t+o(1).
\end{equation}

Suppose that the dimension of $M$ is odd and that the Morse function
is self-indexing (i.e. $f(x)=index(x)$ for any critical point $x$
of $f$). Theorem A of \cite{bfk3} implies that in this  case \refe{bfk}
reduces to the expansion
\eq{as-exp}
      \log\rstt=a_0 + a_1t + b\log \left(\frac t\pi\right) + o(1)
\end{equation}
and also provides formulae for the coefficients $a_0, a_1, b$.

In the present paper we apply the Bismut-Zhang theorem about the
comparison of analytic and combinatorial torsion
(\cite[Theorem 0.2]{bz1}) to prove that
\refe{as-exp} remains true in the general case, when $\dim M$ is not
necessarily odd and the Morse function is not necessarily
self-indexing.

In the case where all the critical values of $f$ are rational
we  calculate explicitly   the coefficients of the asymptotic
expansion \refe{as-exp} in terms of the spectral sequence of a
filtration associated to $f$.  These calculations are based on
the descriptions of the singularities of the torsion obtained
by Farber \cite{fa}.

Now we shall discuss the precise formulations of these results.

\ssec{assum}{Assumptions on $f$, $g^F$ and $g^{TM}$}
Let $f:M\to \RR$ be a Morse function.
Denote by \nf\ the gradient vector
field of $f$ with respect to the metric \gtm. Let $B$ be the
finite set of zeroes of \nf.

We shall assume that the following  conditions are satisfied
(cf. \cite[page 5]{bfk3}):
\begin{enumerate}
    	\item The gradient vector field \nf\ satisfies {\em the Smale
	      transversality conditions} \cite{sm1,sm2} (for any two
	      critical points $x$ and $y$ of $f$ the stable manifold
	      $W^s(x)$ and the unstable manifold $W^u(y)$,
	      with respect to
	      \nf, intersect transversally).
  	\item For any $x\in B$, the metric \gf\ is flat near $B$ and
  	      there is a system of coordinates
	      $y=(y^1\nek y^n)$ centered at $x$ such that near $x$
       {
	\noindent
	    \[
	      	\gtm=\sum_{i=1}^n|dy^i|^2, \quad
	        f(y)=f(x)-\frac12\sum_{i=1}^{index(x)}|y^i|^2+
		\frac12\sum_{i=index(x)+1}^n|y^i|^2.
	     \]
        }
 \end{enumerate}


\ssec{miltor0}{The Milnor torsion} The {\em Thom-Smale complex}
is a finite dimensional complex generated by the fibers of
$F_x \ (x\in B)$ of $F$ whose cohomology is canonically
isomorphic to $H^\b(M,F)$. The metric \gf\ on $F$ defines the
Euclidean structure on the Thom-Smale complex.
The Witten deformation $\gft=e^{-2tf}\gf$ of the metric on
$F$ determines  the deformation of this structure. We shall refer to
this deformation as to the Witten deformation of the Thom-Smale
complex.

The {\it Milnor torsion} \mtt\ is
the torsion of the Thom-Smale complex corresponding to the metric \gft\
on $F$.

\lem{mil-as} The function \mtt\ admits an asymptotic expansion as
   $t\to+\infty$ of the form
   \eq{mil-as}
	\log \mt(t)=\alp+\bet t+o(1).
   \end{equation}
\elem
\refl{mil-as} is proved in \refss{wit-mt}.

\ssec{rational}{The case where all the critical values
of $f$ are rational}
Suppose now  all the critical values $f_1< \cdots<f_l$ of $f$ are
rational. In this situation we shall calculate the coefficients
$\alp, \bet$ of the asymptotic expansion \refe{mil-as}.

Assume that $d\in\NN$ and $p_1\nek p_l\in\ZZ$ are such that
\eq{f=p/d}
	f_i=\frac {p_i}d, \qquad  (1\le i\le l).
\end{equation}
After the change of the deformation parameter
$t\mapsto \tau=e^{-\frac td}$ the Witten deformation
of the Thom-Smale complex satisfies
the conditions of Theorem 6.6 of \cite{fa}.
Hence, the coefficients $\alp, \bet$ of the asymptotic expansion
\refe{mil-as} may be calculated in terms of the spectral sequence
of this deformation (cf. \refss{ssdef}).
This spectral sequence admits the following geometric description.

Let $m,k \in \ZZ$  be such that $(m-k)/d< f(x)< m/d$ for
any $x\in M$. We define the filtration
$\varnothing=U^0\i \dots\i U^k=M$
on $M$ by
\eq{fil}
  	U^i=f^{-1} \left[\frac{m-i}d,\frac md\right]
	  \quad (0\le i\le k).
\end{equation}
In \refs{filtrcom} we show  that the spectral sequence
of the Witten deformation of the Thom-Smale complex my be
expressed in terms of the spectral sequence  $(E^{p,q}_r,d_r)$
associated with this filtration.

\rem{ss}The filtration \refe{fil} and the spectral sequence
$(E^{p,q}_r,d_r)$ depend on the choice of $d$  in \refe{f=p/d}.
Unfortunately, the spectral sequence associated with
seemingly more natural filtration  $\varnothing=
V^0\i \dots\i V^l=M, \ V^i=f^{-1}[f_i,m/d]$
is not connected to the spectral sequence of the Witten deformation.
\erem

Let \sst\ be the torsion of the spectral sequence
$(E^{p,q}_r,d_r)$ (cf. \refd{sst}). Note that our
definition of the torsion of a spectral sequence is
slightly different from \cite{fa} (cf. \refr{sst}).

\prop{mil-as-r}
   If all the critical values of $f$ are rational, then the
   coefficients $\alp$ and $\bet$ in \refe{mil-as} are
   given by the formulae
   \begin{gather}\label{E:mil-as,a}
	\alp=\log \sst;\\  \label{E:mil-as,b}
	\bet=-\frac1d   \bigg( \sum_{p,q\ge 0}
        (-1)^{p+q}{(p+q)}\sum_{r\ge 1}r
	    \big(\dim E^{p,q}_r-\dim E^{p,q}_{r+1}\big)\bigg).
   \end{gather}
\eprop
\refp{mil-as-r} is proved in \refs{mil-as}.

\rem{ration} The assumption that the critical values of $f$ are rational
does not seem  natural. It would be very interesting to obtain formulae
for the coefficients $\alp, \bet$ of \refe{mil-as}  without this
assumption.  These formulae would represent $\alp, \bet$ as functions
of the critical values of the Morse function $f$. Unfortunately, those
functions are not continuous. One can show that they can have jumps when
the numbers $f_1\nek f_l$ are rationally dependent.

Note that, if the critical points of $f$ are not rational,
the substitution $t\mapsto \tau=e^{-\frac td}$ is not defined and,
hence, Farber's theorem can not be applied to the study of
the Witten deformation of the Thom-Smale complex.
\erem

\ssec{notation'}{Notation} Following \cite{bz1}, we introduce the
following definitions.

Let  \ntm\ be the Levi-Civita connection on $TM$ corresponding to the
metric \gtm, and let \etm\ be the associated representative
of the Euler class of $TM$ in Chern-Weil theory.

Let \mqf\ be the Mathai-Quillen (\cite[\S 7]{mq}) $n-1$ current on $TM$
(see also \cite[Section 3]{bgs}\ and \cite[Section IIId]{bz1})
which restriction on $TM/\{0\}$ is induced by a smooth
form on the sphere bundle which transgresses the form $\etm$.

Let $\nabla^F$ denote the flat connection on $F$ and
let \thet\ be the 1-form on $M$ defined by (cf. \cite[Section IVd]{bz1})
\eq{mqform}
  \thet=\Tr\big[(\gf)^{-1}\nabla^F\gf\big].
\end{equation}
Set
\begin{gather}\label{E:chi}
	 	\chi(F)=\sum_{i=0}^n(-1)^i\dim H^i(M,F),\\
	\chi'(F)=\sum_{i=0}^n(-1)^ii\dim H^i(M,F), \label{E:chi'}\\
 	\Tr_s^B[f]=\sum_{x\in B}(-1)^{index(x)}f(x). \label{E:tr(f)}
\end{gather}

Our main result is the following
\th{mainth2}
  Let \rstt\ be the analytic torsion corresponding to the metric
  $\gft=e^{-2tf}\gf$.

  (i) The function $\log\rstt$ admits an asymptotic expansion for
  $t\to +\infty$  of the form
  \eq{mainth}
      \log\rstt=a_0 + a_1t + b\log \left(\frac t\pi\right) + o(1),
  \end{equation}
  where
  \eq{b}
      b=\frac n4 \chi(F)-\frac12\chi'(F)
  \end{equation}
  and $a_0,a_1$ are real numbers depending on $f,\ \gf$ and $\gtm$.

  (ii) Assume that all the critical values of $f$
  are rational and let the
  integer $d$ and the spectral sequence $(E^{p,q}_r,d_r)$ be as in
  \refss{rational}. Then the coefficients $a_0$ and $a_1$ in
  \refe{mainth} are given by the formulae
  \eq{a0}
	a_0=\log\sst -\frac12\int_M\thet(\nf)^*\mqf;
  \end{equation}
  \begin{multline} \label{E:a1}
      a_1=-\rk(F)\int_Mf\etm\\
	 -\frac 1d\sum_{p,q\ge 0}(-1)^{p+q}{(p+q)}\sum_{r\ge 1}r
	    \big(\dim E^{p,q}_r-\dim E^{p,q}_{r+1}\big)+
		\rk(F)\Tr_s^B[f].
  \end{multline}
\eth

\rem{elipt} The Witten deformation is a deformation
of an elliptic complex.  If this deformation were elliptic with
parameter (\cite{sh}, \cite{bfk1}) then, by
\cite[Theorem A.3]{bfk1}, its torsion would have an
asymptotic expansion for $t\to \infty$, whose
coefficients would be given by local expressions.
It fails to be elliptic with parameter precisely at the critical
points of the Morse function $f$. In \cite{bfk3},
Burghelea, Friedlander and Kappeler used the
Mayer-Vietoris formula for elliptic operators (\cite{bfk1})
to show that its torsion continues to have an asymptotic expansion
with  coefficients which are not local anymore. \reft{mainth2}
provides explicit formulae for  these coefficients.
\erem
\rem{dai-mel} The connection between the asymptotic behavior of the
analytic torsion and the spectral sequence associated with the
deformation was discovered by Farber \cite{fa}.
We discuss Farber's results in \refs{1pardef}.

In \cite{dm}, Dai and Melrose have obtained
the asymptotic of the Ray-Singer analytic torsion in the
adiabatic limit. Their result is also expressed in terms of a spectral
sequence.
\erem

\ssec{slfind}{The case where the Morse function is self-indexing}
Suppose now that the Morse function
$f:M\to \RR$ is self-indexing and choose $d=1$ (cf.
\refss{rational}).
For $0\le i\le n$, let $m_i$ denote the number of $x\in B$ of
index $i$. Then  (cf. \cite[page 30]{bz1})  the spectral
sequence $(E^{p,q}_r,d_r)$ degenerates in the second term and
\eq{e1e2}
  \begin{aligned}
	\dim E^{p,q}_1 &=
	     \begin{cases}
		m_p\rk(F), \quad&\text{if}\quad q=0;\\
		0,   \quad&\text{if}\quad q\not=0.
	     \end{cases}  \\
	\dim E^{p,q}_2 &=
	     \begin{cases}
		\dim H^p(M,F), \quad &\text{if} \quad q=0;\\
		0,   \quad&\text{if}\quad q\not=0.
	     \end{cases}
   \end{aligned}
\end{equation}
Hence,
\eq{sind}
	\sum_{p,q\ge 0}(-1)^{p+q}{(p+q)}\sum_{r\ge 1}r
	    \big(\dim E^{p,q}_r-\dim E^{p,q}_{r+1}\big)=
		\rk(F)\Tr_s^B[f]-\chi'(F).
\end{equation}
Also, the torsion \sst\ of the spectral sequence is easily seen to be
equal to the Milnor torsion \mt.

We obtain the following corollary (cf. \cite[Theorem A]{bfk3},
\cite[Theorem 0.4]{b})
\cor{slfind}
  If the Morse function $f:M\to\RR$ is self-indexing, then
  the coefficients $a_0, a_1$ of the asymptotic expansion
  \refe{mainth} are given  by the formulae
  \begin{gather} \label{E:a0'}
	a_0=\log\mt -\frac12\int_M\thet(\nf)^*\mqf;\\   \label{E:a1'}
      a_1=-\rk(F)\int_Mf\etm+\chi'(F).
  \end{gather}
\ecor

\ssec{methpr}{The method of the proof}  Our method is completely
different from that of \cite{bfk3}. In \cite{bfk3}\ the  asymptotic
expansion  is proved by direct analytic arguments and,
then is applied to get a new proof of the
Ray-Singer conjecture \cite{rs} (which was originally proved by
Cheeger \cite{ch} and M\"uller \cite{mu1}).

In \refs{mainres} of the present paper, we use the Bismut-Zhang
extension of this conjecture (\cite[Theorem 0.2]{bz1})
to get the following proposition
\prop{reduce}
   As $t\to +\infty$, the following identity holds
   \begin{multline}\label{E:reduce}
	\log\rstt - \log\mtt=\\
	  -\frac12\int_M\thet(\nf)^*\mqf-t\rk(F)\,\int_Mf\etm \\
		+t\rk(F) \Tr_s^B[f]+
		\left(\frac n4\chi(F)-\frac 12\chi'(F)\right)
                    \log\left(\frac t\pi\right)+ o(1).
   \end{multline}
\eprop

\reft{mainth2} follows from \refp{reduce}, \refl{mil-as} and
\refp{mil-as-r}.

\subsection*{Acknowledgments} I  am very thankful
to Michael Farber for
explaining me the role of  the spectral sequence of a deformation
and for valuable discussions.

\sec{fdcom} {The torsion of a finite dimensional complex}

In this section we follow \cite[Section 1a]{bz1}.
\ssec{detline}{The determinant line} If $\lam$ is a real line, let
$\lam^{-1}$ be the dual line. If $E$ is a finite dimensional
vector space, set
\[
   	\det E= {\bigwedge} ^{\max}(E)
\]
Let
\eq{fd-com1}
	(V^\b,\p):\ 0\to V^0@>\p>>\cdots @>\p>> V^n\to 0
\end{equation}
be a  complex of finite dimensional Euclidean vector spaces. Let
$H^\b(V)=\bigoplus_{i=0}^n H^{i}(V)$ be the cohomology of
$(V^{\b},\p)$.
Set
\begin{gather}\label{E:det}
    \det V^\b=\bigotimes_{i=0}^n \big(\det V^{i}\big)^{(-1)^i}, \\
    \det H^{\b}(V)=\bigotimes_{i=0}^n \Big(\det H^{i}(V)\Big)^{(-1)^i}.
\end{gather}
Then, by \cite{kmu}, there is a canonical
isomorphism of real lines
\eq{kmu-isom}
 	\det H^{\b}(V)\simeq \det V^\b.
\end{equation}

\ssec{2met}{Two metrics on the determinant line}
The Euclidean structure on $V^\b$ defines a metric on $\det V^\b$. Let
$\fdm$ be the metric on the line $\det H^\b(V)$ corresponding to this
metric via the canonical isomorphism \refe{kmu-isom}.

Let $\p^*$ be the adjoint of $\p$ with
respect to the Euclidean structure on \ts. Using the finite dimensional
Hodge theory, we have the canonical identification
\eq{hodge-ts}
 H^i(V^\b,\p)\simeq \{v\in V^i:\ \p v=0, \p^* v=0\},\ \ 0\le i\le n.
\end{equation}
As a vector subspace of $V^i$, the vector space in the right-hand
side of  \refe{hodge-ts} inherits the Euclidean metric. We denote by
\fdmb\ the corresponding metric on $\det H^{\b}(V)$. We shall refer to
this metric as to the {\em Hodge metric} on $\det H^\b(V)$.

The metrics \fdm\ and \fdmb\ do not coincide in general.
We shall describe the discrepancy.

\ssec{fdtor}{The  torsion of a finite dimensional complex}
Set $\Del=\p\p^*+\p^*\p$ and let $\Pi:V^\b\to\Ker \Del$
be the orthogonal projection. Set $\Pi^\perp=1-\Pi$.

Let $N$ and $\tau$ be the operators on $V^\b$ acting on
$V^i\ (0\le i\le n)$ by multiplication by $i$ and $(-1)^i$ respectively.
If $A\in \End(V^\b)$, we define the supertrace $\Tr_s[A]$ by the formula
\eq{suptr}
	\Tr_s[A]=\Tr[\tau A].
\end{equation}
{}For $s\in\CC$, set
\[
    \zet^V(s)=-\Tr_s\big[N(\Del)^{-s}\Pi^\perp\big].
\]

\defe{fdtor}
  The {\em  torsion} of the complex
$(V^\b,d)$ is the number
  \eq{fdtor}
	\rho=\exp\Big(\frac12\frac{d\zet^V(0)}{ds}\Big).
  \end{equation}
\edefe
We denote by $\Del^i\ (0\le i\le n)$ the restriction of $\Del$ on $V^i$ .
Let $\{\lam^i_j\}$ be the set of nonzero eigenvalues of $\Del^i$. Then
\eq{tor=det}
	\log \rho = \frac12\sum_{i,j} (-1)^ii\log \lam^i_j.
\end{equation}

The following result is proved in \cite[Proposition 1.5]{bgs}
\eq{fdm-fdmb}
	\fdm=\fdmb \cdot\rho.
\end{equation}

\sec{milmet}{The Milnor metric and the Milnor torsion}

In this section we recall the definitions of the Milnor metric and
the Milnor torsion and prove \refl{mil-as}.
\ssec{mildetl}{The determinant line of the cohomology}
Let $\h{\b}= \bigoplus_{i=0}^n\h{i}$ be the cohomology
of $M$ with coefficients in $F$ and
let \dh{\b} be the line
\eq{mildetl}
	\dh{\b}=\bigotimes_{i=0}^n \Big(\dh{i}\Big)^{(-1)^i}.
\end{equation}

\ssec{th-sm}{The Thom-Smale complex}
Let $f:M\to\RR$ be a Morse function satisfying the
Smale transversality conditions \cite{sm1,sm2}\ (for any two
critical points $x$ and $y$ of $f$ the stable manifold
$W^s(x)$ and the unstable manifold $W^u(y)$, with respect to
\nf, intersect transversally).

Let $B$ be the set of critical points of $f$.   If $x\in B$,
 let $F_x$ denote the fiber of $F$ over $x$ and let $[W^u(x)]$
denote the real line generated by $W^u(x)$. For $0\le i\le n$,
set
\eq{th-sm}
  C^i(W^u,F)=
    \bigoplus_{\begin{Sb} x\in B\\ index(x)=i \end{Sb}}
		[W^u(x)]^*\ten_{\RR} F_x.
\end{equation}
By a basic result of Thom (\cite{thom}) and Smale (\cite{sm2})
(see also \cite[pages 28--30]{bz1}), there are well defined
linear operators
\[ \p:C^i(W^u,F)\to C^{i+1}(W^u,F), \]
such that  the pair $(C^\b(W^u,F),\p)$ is a complex  and there is
a canonical identification of $\ZZ$-graded vector spaces
\eq{tsc=c}
 	  H^\b(C^\b(W^u,F),\p)\simeq H^\b(M,F).
\end{equation}
\ssec{milmet}{The Milnor metric}
By \refe{kmu-isom} and \refe{tsc=c}, we know that
\eq{km-ts}
	\dh{\b}\simeq \dts.
\end{equation}
The metric \gf\ on $F$ determines the structure of an Euclidean vector
space on \ts. This structure induces a metric on $\det\ts$.
\defe{milmet}
   The {\em Milnor metric} \mm\ on the line \dh{\b}\ (cf.
   \cite[Section Id]{bz1}) is the metric
   corresponding to the above metric on \dts\ via the
   canonical isomorphism \refe{km-ts}.
\edefe
\rem{milinv}By Milnor \cite[Theorem 9.3]{mi1}, if \gf\ is a flat
  metric on $F$, then
  the Milnor metric  coincides with the Reidemeister metric defined
  through a smooth triangulation of $M$. In this case \mm\
  does not depend upon $f$ and \gtm\ and, hence, is a
  topological invariant of the flat
  Euclidean vector bundle $F$.
\erem

\defe{mmb} The {\em Hodge-Milnor metric} \mmb\ on \dh{\b}\
is the metric corresponding to the Hodge metric  (cf. \refss{2met}) on
$\det H^\b (\ts,\p)$ via \refe{km-ts}.
\edefe

\defe{miltor}
The {\em Milnor torsion} is the  torsion of the Thom-Smale complex
(cf. \refd{fdtor}).
\edefe
{}From \refe{fdm-fdmb}, we obtain
\eq{mm-mmb}
	\mm=\mmb \cdot\mt.
\end{equation}

\ssec{wit-mil}{Deformation of the Milnor metric}
The metric \mm\ depends on the metric \gf.  Let $\gft=e^{-2tf}\gf$
and let \mmt\ be the corresponding Milnor metric.
Recall from \refss{notation'} the notation
\eq{tr}
\Tr^B_s[f]=\sum_{x\in B}(-1)^{index(x)}f(x).
\end{equation}
Obviously,
\eq{wit-mil}
	\mmt= e^{-t\rk(F) Tr_s^B[f]}\cdot\mm.
\end{equation}

\ssec{wit-mt}{Proof of \refl{mil-as}}
Assume that  $g_t$ is the metric on \ts\ induced by the metric
$\gft=e^{-2tf}\gf$ on $F$. Let $\calF\in\End(\ts)$ which,
for  $x\in B$, acts on $[W^u(x)]^*\otimes F_x$ by multiplication
by $f(x)$. Then, for any $x,y\in\ts$ , we have
\eq{gt}
	g_t(x,y)=g_0(e^{-2t\calF}x,y).
\end{equation}

Let $\p^*_t$ be the adjoint of $\p$ with respect to the metric
$g_t$.  Clearly,
\eq{p-pt}
	\p_t^* =e^{2t\calF}\p^*_0e^{-2t\calF}.
\end{equation}
Set $\Del_t=\p\p^*_t+\p^*_t\p$ and denote by $\Del_t^i$  the
restriction of $\Del_t$ on $C^i(W^u,F)$. The number
$k^i= \dim\Ker\Del^i_t$ does not depend on $t$.
Hence, the characteristic polynomial $\det(\Del_t^i-xI)$
of $\Del^i_t$ may be written in the form
\eq{chpoly}
	\det(\Del^i_t -xI)=x^{k^i}
		\sum_{j=0}^{\dim V^i-k^i}
			a_j^i(t)x^j,
\end{equation}
where $a^i_0(t)\not=0$ is equal to the product of the nonzero
eigenvalues  of $\Del^i_t$.

Let \mtt\ denote the Milnor torsion corresponding to the metric \gft.
{}From \refe{tor=det}, we see that
\eq{tor=alp}
	\log\mtt=\frac12\sum_{i=0}^n (-1)^iia^i_0(t).
\end{equation}
By \refe{p-pt}, $a^i_0(t)$ is a polynomial in $e^{-2tf_1}\nek
e^{-2tf_l}$. Hence, for any  $0\le i \le n$, there exist
real numbers $\alp^i, \bet^i$  such that, for $t\to +\infty$,
\eq{a=alpt}
	a^i_0(t)=\alp^i e^{\bet^i t}(1+o(1)).
\end{equation}
{}From  \refe{tor=alp} and \refe{a=alpt}, we obtain
\eq{mil-as''}
	\log\mtt=\frac12\sum_{i=0}^n (-1)^ii\log \alp^i+
		\frac12 t\sum_{i=0}^n (-1)^ii\bet^i+o(1),
\end{equation}
completing the proof of the lemma.

\sec{ssdef}{The spectral sequence of a deformation}

\ssec{def}{A deformation of a complex}
Let
\eq{com}
	(V^\b,\p_0):\  0\to V^0@>\p_0>>\cdots@>\p_0>>V^n\to 0
\end{equation}
be a complex of real vector spaces.

Denote by $\calO$ the ring of germs at
the origin of real analytic functions of one  variable $t$. By
a {\em deformation of $(V^\b,\p_0)$} we shall understand a complex
\eq{calcom}
	(\calV^\b,\p):\  0\to \calV^0@>\p>>\cdots@>\p>>\calV^n\to 0
\end{equation}
of $\calO$-modules together with a fixed isomorphism between the fiber
$\calV^\b\otimes \calO/t\calO$ of $(\calV^\b,\p)$
at the point $t=0$ and $(V^\b,\p_0)$.

\ssec{ssdef}{The spectral sequence} Given a deformation \refe{calcom}
of a complex \refe{com}, we consider a short exact sequence of complexes
\eq{vvb}
	0\to\calV^\b@>t>>\calV^\b@>q>>\calB^\b\to 0.
\end{equation}
Here $t:\calV^\b\to\calV^\b$ is the multiplication by $t$,
$\calB^\b$ is the quotient complex $\calB^\b= \calV^\b/t\calV^\b$
and $q$ is the quotient map.

The exact sequence \refe{vvb} induces a long exact sequence of
cohomology
\eq{calcoh}
	\cdots\to H^m(\calV^\b)@>t>>H^m(\calV^\b)@>q>>
	   H^m(\calB^\b)@>r>>H^{m+1}(\calV^\b)\to\cdots,
\end{equation}
which may be rewritten as an exact couple
\eq{calex2}
    \begin{CD}
	H^\b(\calV^\b)&&\longrightarrow&&H^\b(\calV^\b)\\
	&\nwarrow && \swarrow &\\
	&&H^\b(\calB^\b)&&
    \end{CD}
\end{equation}
According to the standard rules the latter generates a (Bockstein)
spectral sequence $(\calE^i_r,d_r)$. To describe it explicitly
set
\eq{z}
	\calZ^i_r=\big\{s\in\calV^i:\ \p s\in t^r\calV^{i+1} \big\}.
\end{equation}

Then
\begin{gather}
  	\calE^i_r=\left\{
       	   \begin{aligned}
	      &V^i \quad \text{for} \quad r=0, \\
	      &\calZ^i_r/(t\calZ^i_{r-1}+
		   t^{1-r}\p \calZ^{i-1}_{r-1})\quad
	      \text{for} \quad r> 0,
	    \end{aligned}
	\right.\\
\intertext{and the differential}
	d_r:\calE^i_r\to \calE^{i+1}_r
\end{gather}
is the homomorphism induced by the action of $t^{-r}\p$
on $\calZ^i_r$. Then $H^\b(\calE_r^\b,d_r)\simeq
\calE^\b_{r+1}$, i.e. one gets a
spectral sequence.

Note that the sequence $(\calE^i_r,d_r)$ is completely determined by
\refe{vvb}.

\sec{1pardef}{Deformation of the torsion of
	a finite dimensional complex}

In this section we consider  a one parameter family $(V,\p_t)$ of
finite dimensional Euclidean complexes. This family can be
considered as a deformation of a complex and, according to
the previous section, gives rise to a spectral sequence.
We define the torsion of this spectral sequence. Finally we prove
a theorem by Farber
(cf. \cite[Theorem 6.6]{fa}) which describes the asymptotic for
$t\to 0$  of the torsion $\rho(t)$ of the complex $(V,\p_t)$.

In this section we essentially follow \cite{fa}.

\ssec{1pardef}{A family of Euclidean complexes} Let
\eq{1pardef}
	(V^\b,\p_t):\  0\to V^0@>\p_t>>\cdots@>\p_t>>V^n\to 0
\end{equation}
be a one parameter family of  complexes of  finite
dimensional Euclidean vector spaces. We shall assume that the operators
\[
    \p_t:V^i\to V^i
\]
depend analytically on a parameter $t$ varying within an interval
$(-\eps,\eps)$. That means  that $\p_t$ may be represented as a
convergent power series
\eq{series}
	\p_t=\p_0+t\p_1+\cdots
\end{equation}
with coefficients in $\End(V)$.

\ssec{germ}{The germ complex}
Set $\calV^i=
\calO\ten_{\RR} V^i\ (0\le i\le n)$. The family \refe{1pardef}
can be understood as a single complex
\eq{germcom}
	(\calV^\b,\p):\  0\to \calV^0@>\p>>\cdots@>\p>>\calV^n\to 0
\end{equation}
of $\calO$-modules, where the differential $\p:\calV^\b\to\calV^\b$
is given by
\eq{tild}
	\p s(t)=\p_ts(t)\quad (s(t)\in \calV).
\end{equation}
The complex \refe{germcom} is called  the {\em germ complex}
(cf. \cite[Section 2.5]{fa}).

\ssec{kato}{The parameterized spectral decomposition}
Let $\p_t^*\ (-\eps<t<\eps)$
be the adjoint of $\p_t$ with respect to the Euclidean structure on
$V^\b$. Set $\Del_t=\p_t\p_t^*+\p_t^*\p_t$. For $0\le i\le n$, we denote
by $\Del^i_t$ the restriction of $\Del_t$ on $V^i$.

By a theorem of Rellich \cite[\S 1, Theorem 1]{re} (see
also \cite[Ch. 7, Theorem 3.9]{kato}), there exists a family
of analytic curves $\phi^i_j(t)\in\calV^i\ (1\le j\le
\dim V^i)$ and a sequence of real valued analytic functions
$\lam^i_j(t)\in\calO\ (1\le j\le \dim V^i)$ such that for any value
of $t$ the numbers $\{\lam^i_j(t)\}$ represent all the repeated
eigenvalues of $\Del^i_t$ and $\{\phi^i_j(t)\}$ form a complete
orthonormal basis of corresponding eigenvectors of $\Del^i_t$.

As the operators $\Del_t^i$ are non negative for any $t$, the
functions $\lam^i_j(t)$ depend only on $t^2$.

Suppose that $\lam^i_j(t)$  and $\phi^i_j(t)$  have been numerate
so that there exist integers $0=N^i_0\le N_1^i\le
\dots \le N^i_{m_i}\le N^i_{m_i+1} = \dim V^i$ such that
\begin{enumerate}
    %
    \item  $\lam^i_j(t)=t^{2k} \olam^i_j(t)$ with $\olam^i_j(0)\not=0$
		for $N^i_k+1\le j\le N^i_{k+1,}\ 0\le k\le m_i-1 $;
    \item  $\lam^i_j(t)\equiv 0$ for $j\ge N^i_{m_i}+1$.
\end{enumerate}

\ssec{tor(t)}{The torsion as a function of the parameter}
{}For each $t\in (-\eps,\eps)$, we shall denote by $\rho(t)$
the torsion of the complex $(V^\b,\p_t)$. The following lemma
follows directly from \refe{tor=det}.
\lem{tor(t)}
   The function $\rho(t)$ admits an asymptotic expansion for $t\to 0$
   of the form
   \eq{tor(t)}
	\log\rho(t)=  \frac 12
	  \sum_{i=0}^n(-1)^ii\sum_{j=1}^{N^i_{m_i}}\log\olam^i_j(0)+
	    \bigg(
	     \sum_{i=1}^n (-1)^ii\sum_{k=1}^{m_i-1} k(N_{k+1}^i-N_k^i)
	     \bigg) \log(t)+o(1).
   \end{equation}
\elem

Now our goal is to express the right hand side of \refe{tor(t)} in
terms of the spectral sequence of deformation \refe{germcom}.

\ssec{dphi}{}
Let $\calV^i_k \ ( 0\le k\le m_i)$ denote the
submodule of $\calV^i$ generated by the set $\big\{\phi^i_j|\,
N_k^i+1\le j\le N^i_{k+1}\big\}$. Since the
operator $\p_t$ commutes with
$\Del_t$ for any $t$, we get $\p\calV^i_k\i \calV^{i+1}_k$.
Then the equality
\eq{del-dd}
	\langle\Del_ts(t),s(t)\rangle =
		\|\p_t s(t)\|^2+\|\p_t^*s(t)\|^2
\end{equation}
implies the following lemma.
\lem{dphi}If $N_k^i+1\le j\le N^i_{k+1}$, then
   \eq{dphi}
	\p\phi^i_j\in t^k\calV^{i+1}_k.
   \end{equation}
\elem


\ssec{hodge}{The Hodge spectral sequence} Now we shall give a Hodge
theoretical description of the  spectral sequence $(\calE^i_r,d_r)$
associated with the deformation \refe{germcom}.

{}For $r\ge 0$ and $ 0\le i \le n$, set
\eq{h i r}
	H^i_r=\SPAN\big\{\phi^i_j(0):\
		N_r^i+1\le j\le \dim V^i\big\}\i V^i.
\end{equation}

By \refl{dphi},
$\sum_{j=N^i_r+1}^{\dim V^i} a^i_j \phi^i_j(t)\in \calZ^i_r$ for
any numbers $a^i_j\in \RR$. Hence, there is a natural function
 $\Phi_r:H^\b_r\to \calE^\b_r$ which maps
$\sum_{j=N^i_r+1}^{\dim V^i} a^i_j \phi^i_j(0)$ to the image
of $\sum_{j=N^i_r+1}^{\dim V^i} a^i_j \phi^i_j(t)\in \calZ^i_r$ in
$\calE^i_r= \calZ^i_r/(t\calZ^i_r+t^{1-r}\p \calZ^{i-1}_{r-1})$.
\lem{injec}
  {}For any $r\ge 0$, the map $\Phi_r:H^\b_r\to \calE^\b_r$
  is injective.
\elem
\prf
  Suppose that
  \[
	\Phi_r\bigg(\sum_{j=N^i_r+1}^{\dim V^i}
		a^i_j \phi^i_j(0)\bigg)=0.
  \]
  Then
  \[
	\sum_{j=N^i_r+1}^{\dim V^i} a^i_j \phi^i_j(t)=
		t\alp(t)+t^{1-r}\p\bet(t),
  \]
  where $\alp\in \calZ^i_r,\ \bet\in \calZ^{i-1}_{r-1}$. Hence,
  \eq{injec1}
	\sum_{j=N^i_r+1}^{\dim V^i} a^i_j \phi^i_j(0)=
		t^{1-r}\p\bet(t)\Big|_{t=0}.
  \end{equation}
  By \refl{dphi},
  \eq{injec2}
	t^{1-r}\p \bet(t)\Big|_{t=0}\in
		\SPAN\big\{\phi^i_j(0):\
			N^i_{r-1}+1\le j\le N^i_r \big\}.
  \end{equation}
  Hence, \refe{injec1} implies that
  $\sum_{j=N^i_r+1}^{\dim V^i} a^i_j \phi^i_j(0)= 0$.
\epr

By \refl{dphi}, there is a map $\del_r:H^\b_r\to H^\b_r$
defined by
\eq{del}
	\del_r:\sum_{j=N^i_r+1}^{\dim V^i} a^i_j \phi^i_j(0)=
	t^{-r}\p
	\sum_{j=N^i_r+1}^{\dim V^i} a^i_j \phi^i_j(t)
	\Big|_{t=0}.
\end{equation}
Clearly, $\del_r^2=0$.

Let $\del_r^*:H^\b_r\to H^\b_r$ be the adjoint of $\del_r$.
Using \refe{del} and \refl{dphi}, we see that
\eq{Delphi}
	(\del\del^*+\del^*\del)\phi^i_j(0)=
	   \begin{cases}
		0, \quad &\text{if}\quad j > N^i_{r+1};\\
		\olam^i_j(0), &\text{if}\quad N^i_r+1\le j\le N^i_{r+1}.
	   \end{cases}
\end{equation}
{}From \refe{h i r}, \refe{Delphi}, we get
\eq{hr-hodge}
	H^\b_{r+1}=\big\{h\in H^\b_r:
		\ (\del\del^*+\del^*\del)h=0\big\}=
		      \big\{h\in H^\b_r:\ \del_rh=\del_r^*h=0\big\}.
\end{equation}

The following proposition is equivalent to Theorem 3.3 of \cite{kk}.
A particular case of this result was proved by Forman
(\cite[Theorem 6]{fo}).
\prop{kk} For all $r\ge 0$,

  (i) The map $\Phi_r:H^\b_r\to \calE^\b_r$ is an isomorphism.

  (ii) $\del_r=\Phi_r^{-1}d_r\Phi_r$.

  \eprop
\prf
  {}Following \cite{kk}, we shall prove the proposition by induction
  on $r$.

  The case $r=0$ is obvious. For the inductive step, assume we
  have proven that $\Phi_r:H^\b_r\to \calE^\b_r$ is an isomorphism and
  $\del_r= \Phi_r^{-1}d_r\Phi_r$. Then, by \refe{hr-hodge} and
  the finite dimensional
  Hodge theory, $\calE^\b_{r+1}$ is isomorphic to
  \[
  	\big\{h\in H^\b_r:\ \del_rh=\del_r^*h=0\big\}= H^\b_{r+1}.
  \]
  Then \refl{injec} implies that
  $\Phi_{r+1}:H^\b_{r+1}\to \calE^\b_{r+1}$
  is an isomorphism. From the definition of $d_r$ and $\del_r$,
  we obtain   $\Phi_{r+1}\del_{r+1}=d_{r+1}\Phi_{r+1}$, completing
  the induction step and the theorem.
\epr
\cor{dimE}  For any $0\le i\le n$ and $r\ge 0$,
  \eq{dimE}
	\dim \calE^i_r=\dim V^i-N^i_r.
  \end{equation}
\ecor
As another corollary we obtain the following result by Farber
\cite[Theorem 1.6]{fa}
\prop{spseq}
   The spectral sequence $(\calE^\b_r,d_r)$ stabilizes and the
   limit term $\calE^\b_{\infty}$
  is isomorphic to the cohomology $H^\b(V^\b,\p_t)$ for a generic point
  $t$.
\eprop


\ssec{sst}{The torsion of the spectral sequence}
As a subspace of $V^\b$
the vector space $H^\b_r$ $(r\ge 0)$ inherits the Euclidean metric.
We denote by $\rho_r$ the torsion of the complex $(H^\b_r,\del_r)$
corresponding to this metric. Note that, by \refp{spseq},
$\rho_N=1$ for sufficiently large $N$.

\defe{sst} The {\em torsion of spectral sequence} $(\calE^\b_r,d_r)$
  is the product
  \eq{sst}
  	\sst=\rho_0\rho_1\dots\rho_N,
  \end{equation}
  where $N$ is a sufficiently large number.
\edefe
\rem{sst} The torsion of the spectral sequence of a deformation
  was defined by Farber
  \cite[Section 6.5]{fa} in slightly different terms.
  Note that the torsion of a spectral sequence,
  as it is defined in \cite{fa}, corresponds, in our terms, to the
  product $\rho_1\rho_2\dots\rho_N$.
\erem


\ssec{farber}{Farber theorem}
Using \cite[Proposition 6.3]{fa}, one can easily see that
the following theorem is equivalent to \cite[Theorem 6.6]{fa}.

%
\th{farber}
  The function $\rho(t)$ admits an asymptotic expansion
  for $t\to 0$ of the form
  \eq{farber}
	\log \rho(t)=\log \sst+ \bigg( \sum_{i=0}^n (-1)^ii
		\sum_{r\ge 1}r\big(\dim \calE^i_r-\dim
		   \calE^i_{r+1}\big)\bigg)\log(t)+ o(1).
  \end{equation}
\eth
\prf
  By \refe{tor=det}, \refe{Delphi}, we obtain
  \eq{rho=lam}
	\log \rho_r = \frac 12\sum_{i=0}^n(-1)^ii
	   \sum_{i={N_r^i+1}}^{N_{r+1}^i} \log\olam^i_j(0).
  \end{equation}
  From \refe{tor(t)}, \refe{dimE}, \refe{sst}  and \refe{rho=lam},
  we get \refe{farber}.
\epr

\sec{filtrcom}{Deformation of a filtered complex}

In this section we apply \reft{farber} to a deformation of a
filtered complex. The results of this section are closely related to
\cite[Theorem 9]{fo}.
\ssec{ssfil}{The spectral sequence of a filtration}
Suppose that
\eq{fd-com1'}
	(V^\b,\p):\ 0\to V^0@>\p>>\cdots @>\p>> V^n\to 0
\end{equation}
is a filtered complex of real  vector spaces
with an increasing filtration
\eq{filtr'}
	0=F_0V^\b\i F_1V^\b\i\dots, \qquad V^\b=\bigcup_{i\ge 0}F_iV^\b.
\end{equation}
Denote by $(E^{p,q}_r,d_r)$ $(r\ge 0)$ the spectral
sequence of this filtration and set $E_r^i=
\bigoplus_{p+q=i}E^{p,q}_r$.
Then   $(E^i_r,d_r)$ is also a spectral sequence. We shall need
the following description of $(E^i_r,d_r)$ (cf. \cite[\S 14]{botu}).

Let $i:\bigoplus_{p=0}^\infty F_pV^\b\to
\bigoplus_{p=0}^\infty F_pV^\b$  be
the map induced by the inclusions $F_pV^\b\hookrightarrow
{}F_{p+1}V^\b$ $(p\ge 0)$. Consider the exact sequence of complexes
\eq{ffb}
	0\to\bigoplus_{p=0}^\infty  F_pV^\b @>i>>
	     \bigoplus_{p=0}^\infty F_pV^\b @>j>>B^\b\to 0.
\end{equation}
The complex $B^\b$ is isomorphic to the associated graded complex
$grV^\b$ of $V^\b$.

In a standard way, the exact sequence \refe{ffb} gives rise to
a spectral sequence, which is isomorphic to $(E^i_r,d_r)$.
It follows, that $(E^i_r,d_r)$ is completely determined by
\refe{ffb}.

\ssec{rees}{The Rees complex}We shall construct a complex
$(\calV^\b,\overline{\p})$ of $\calO$-modules as follows
\eq{rees}
  \begin{aligned}
	\calV^i=
	   \big\{\sum_{m=0}^N v_mt^m :\  v_i\in &F_mV^i, N\in\NN\big\},
		\qquad \calV^\b=\bigoplus_{i=0}^n \calV^i,\\
	\overline{\p}\sum_{m=0}^N v_pt^p &=\sum_{m=0}^N (\p v_p)t^p.
  \end{aligned}
\end{equation}
Note that the fiber of $\calV^\b$ at the point $t=0$ is isomorphic
to the associated graded complex $gr V^\b$ of $V^\b$. Hence,
$(\calV^\b,\overline{\p})$ is a deformation of $gr V^\b$ in
the sense of \refss{def}.

As in \refss{ssdef}, we construct a short exact sequence of complexes
\eq{vvb'}
	0\to\calV^\b@>t>>\calV^\b@>q>>\calB^\b\to 0,
\end{equation}
which induces a spectral sequence $(\calE^i_r,\od_r)$.

Note that, as a real vector space, $\calV^\b$ is isomorphic to
$\bigoplus F_kV^\b$ and the multiplication by $t$ corresponds
under this isomorphism to the inclusion $i:\bigoplus F_kV^\b
\hookrightarrow \bigoplus F_kV^\b$. Hence, the exact sequences
\refe{ffb} and \refe{vvb'} are isomorphic. Then so are the spectral
sequences $(E^i_r,d_r)$ and $(\calE_r^i,\od_r)$.


\ssec{met}{Deformation of the torsion}
Suppose that the dimensions of the spaces $V^1\nek  V^n$
are finite. Then there exists $k\in \NN$ such that
$F_jV^\b=V^\b$ for any $j\ge k$. Let $g^{V^0}\nek g^{V^n}$ be
Euclidean metrics on $V^0\nek   V^n$.
With these assumptions we shall present an equivalent description
of deformation \refe{rees}.

Equip $V=\bigoplus_{i=0}^nV^i$ with the metric
$g^V=\bigoplus_{i=0}^n g^{V^i}$, which is the orthogonal
sum of the metrics $g^{V^0}\nek g^{V^n}$.

Let $\Pi_j:V^\b\to F_jV^\b\ (0\le j\le k)$ be the orthogonal
projection. For $t>0$, set
\eq{A}
	A_t=\sum_{j=1}^k t^j(\Pi_{j}-\Pi_{j-1}).
\end{equation}
The operator $A_t$ is invertible for any $t\not= 0$.

Define
\eq{pt}
	\p_t=
	 \begin{cases}
	      A_t^{-1}\p A_t, \quad &\text{for}\quad t\not=0;\\
	      \sum_{j=1}^k(\Pi_j-\Pi_{j-1})\p (\Pi_j-\Pi_{j-1}),
	         \quad &\text{for}\quad t=0.
	  \end{cases}
\end{equation}
As in \refss{1pardef}, the family of complexes $(V_\b,\p_t)$ may be
considered as a single complex $(\tilV,\tp)$ of $\calO$-modules.
Denote by $A:\tilV^\b\to \calV^\b$ the map defined by the formula
\eq{tilv-calv}
	A\sum_{i=0}^Nx_it^i= \sum_{i=0}^N(Ax_i)t^i.
\end{equation}
Then $A$ is an isomorphism of complexes of $\calO$-modules.
Hence, the spectral
sequence of deformation $(\tilV^\b,\tp)$ is isomorphic
to $(E^i_r,d_r)$. Let $\sst$ denote the torsion of this
spectral sequence (cf. \refd{sst}). From \reft{farber}, we get

\prop{filcom}
  Let $\rho(t)$ be the torsion of the complex $(V^\b,\p_t)$ associated
  to the metric $g^V$. Then $\rho(t)$ admits an asymptotic expansion
  for $t\to 0$ of the form
  \eq{fa}
	\log \rho(t)=\log \sst+
          \bigg( \sum_{p,q\ge 0} (-1)^{p+q}{(p+q)}\sum_{r\ge 1}r
	    \big(\dim E^{p,q}_r-\dim E^{p,q}_{r+1}\big)\bigg)
		\log(t) +o(1).
  \end{equation}
\eprop


\ssec{defmet}{Deformation of the metric} Fix $m\in\ZZ$ and
let $g^V_t$ $(t>0)$ be the metric on $V$ defined by the formula
\eq{gVt}
	g_t^V(x,y)= g^V(t^{2m}A_t^{-2} x,y).
\end{equation}

We denote by $\p^*_t$  the adjoint of \p\ with respect to the metric
$g^V_t$. Then $\p^*_t=A_t^2\p^* A^{-2}_t$.  Set
\eq{Del}
	\Del_t=\p\p_t^*+\p_t^*\p.
\end{equation}

Let $\tilpat_t^*$ be the adjoint of $\p_t$ with
respect to the metric $g^V$. Then $\tp^*_t=A_t\p^* A_t^{-1}$. Set
\eq{tilDel}
	\tilDel_t=\p_t\tilpat^*_t+\tilpat^*_t\p_t.
\end{equation}
It is easy to see that
	$\tilDel_t=A_t^{-1}\Del_t A_t.$
Hence, by the definition of the torsion, we obtain the following lemma.
\lem{td-d}
  {}For any $t>0$, the torsion of the complex $(V^\b,\p)$
   associated to the
   metric $g^V_t$ is equal to the torsion  of the
   complex $(V^\b,\p_t)$ associated  to the metric $g^V$.
\end{lemma}

\sec{mil-as}{Proof of \refp{mil-as-r}}

In this section  we assume that all the  critical values of $f$ are
rational.

\ssec{fil-ts}{A filtration on the Thom-Smale complex} Recall
that the integers $d,m$ and $k$ were defined in
\refss{rational}. The Thom-Smale
complex $(\ts,\p)$ possesses a natural filtration
\eq{filt-ts}
	0=F_0\ts\i F_1\ts\i\dots\i  F_k\ts=\ts,
\end{equation}
with
\eq{Fts}
	F_i\ts=\bigoplus_{
		\begin{Sb}
			x\in B\\
			f(x)\ge \frac{m-i}d
		\end{Sb}}
		   [W^u(x)]^*\ten_{\RR} F_x   \qquad (0\le i\le k).
\end{equation}
We denote by $(E^{p,q}_r,d_r)$ the spectral sequence of this filtration.
This spectral sequence is isomorphic to
the spectral sequence of  filtration \refe{fil}.

Let $\sst$ denote the torsion of the spectral sequence
 $(E^{p,q}_r,d_r)$.

\ssec{mil-as-r}{Proof of \refp{mil-as-r}} Recall that $\calF\in
\End(\ts)$ was defined in \refss{wit-mt}. Let
$\Pi_j:\ts\to F_j(\ts)\ (0\le j\le k)$ be the orthogonal
projection.
Set $\tau=e^{-\frac td}$. Then $\tau\to 0$ as $t\to +\infty$ and
\eq{at-F}
	e^{-2t\calF}=
		\tau^{2m}\bigg(\sum_{j=1}^k\tau^j
	   	  (\Pi_j-\Pi_{j-1})\bigg)^{-2}.
\end{equation}
Hence, by \refp{filcom}, \refl{td-d} and  \refe{gt}, we see that,
for $t\to+\infty$,
\eq{def-mt}
	\log \mtt=\log \sst-
          \bigg( \sum_{p,q\ge 0} (-1)^{p+q}{(p+q)}\sum_{r\ge 1}r
	    \big(\dim E^{p,q}_r-\dim E^{p,q}_{r+1}\big)\bigg)
		\frac td+o(1),
\end{equation}
which is exactly \refp{mil-as-r}.

\sec{rsmet}{The Ray-Singer metric and the Ray-Singer torsion}

\ssec{badrsmet}{The $L_2$ metric on the determinant line}
Let $(\om,d^F)$ be the de Rham complex of the smooth sections of
$\bigwedge(T^*M)\ten F$ equipped with the coboundary operator $d^F$. The
cohomology of this complex is canonically isomorphic to \hb.

Let $\ast$ be the Hodge operator associated to the metric \gtm.
We equip \om\ with the inner product
\eq{dr-inpr}
	\langle\alp,\alp'\rangle_{\om}=
		\int_M\langle\alp\land\ast \alp'\rangle_{\gf}.
\end{equation}
By Hodge theory, we can identify \hb\ with the space of harmonic forms
in \om. This space inherits the Euclidean product \refe{dr-inpr}.
The {\em $L_2$ metric} \rsmb\ on \dhb\ is the metric induced
by this product.

\ssec{rstor}{The Ray-Singer torsion}
Let $d^{F*}$ be the formal adjoint of $d^F$ with respect
to the metrics \gtm\ and \gf.

Set $\Del=d^Fd^{F\ast}+d^{F\ast}d^F$ and let $P:\om\to\Ker \Del$ be
the orthogonal projection. Set $P^\perp=1-P$.

Let $N$ be the operator defining the $\ZZ$-grading of \om, i.e. $N$ acts
on $\Ome^i(M,F)$ by multiplication by $i$.

If an operator $A:\om\to \om$ is trace class,
we define its supertrace $\Tr_s[A]$ as in \refe{suptr}.

{}For $s\in\CC$, $\RE\, s >n/2$, set
\[  \zet^{RS}(s)=-\Tr_s\big[N(\Del)^{-s}P^\perp\big].  \]

By a result of Seeley \cite{se}, $\zet^{RS}(s)$ extends to a
meromorphic function of $s\in\CC$, which is holomorphic at $s=0$.

\defe{rstor}
The {\em Ray-Singer torsion \/} is the number
\eq{rstor}
	\rst=\exp\left(\frac12\frac{d\zet^{RS}(0)}{ds}\right).
\end{equation}
\edefe

\ssec{rsmet}{The Ray-Singer metric}
We now remind the following definition
(cf. \cite[Definition 2.2]{bz1}):
\defe{rsmet}
  The {\em Ray-Singer metric \/} \rsm\ on the line \dhb\ is the product
\eq{rsmet}
  \rsm=\rsmb\cdot\rst.
\end{equation}
\edefe

\rem{}When $M$ is odd dimensional, Ray and Singer
  \cite[Theorem 2.1]{rs}\
  proved that the metric \rsm\ is a topological invariant, i.e. does not
  depend on the metrics \gtm\ or \gf.
  Bismut and Zhang \cite[Theorem 0.1]{bz1}\
  described explicitly the dependents of \rsm\ on \gtm\ and \gf\ in
  the  case when $\dim M$ is even.
\erem

\ssec{mm=rm}{Bismut-Zhang theorem}
%
%
%
%
Let  \ntm\ be the Levi-Civita connection on $TM$ corresponding to the
metric \gtm, and let \etm\ be the associated representative
of the Euler class of $TM$ in Chern-Weil theory.

Let \mqf\ be the Mathai-Quillen (\cite[\S 7]{mq}) $n-1$ current on $TM$
(see also \cite[Section3]{bgs}\ and \cite[Section IIId]{bz1}).

Let $\nabla^F$ be the flat connection on $F$ and
let \thet\ be the 1-form on $M$ defined by (cf. \cite[Section IVd]{bz1})
\eq{mqform'}
  \thet=\Tr\big[(\gf)^{-1}\nabla^F\gf\big].
\end{equation}

Now we remind the  following theorem by Bismut and Zhang
\cite[Theorem 0.2]{bz1}.
\begin{Thm}[Bismut-Zhang]\label{T:bz}
The following identity holds
   \eq{bz}
	  \log\left(\frac\rsm\mm\right)^2=
		-\int_M\thet(\nf)^*\mqf.
   \end{equation}
\end{Thm}

\ssec{dfrm-bz}{Dependence on the metric}
The metrics \rsm\ and \mm\ depend,
in general, on the metric \gf. Let $\gft=e^{-2tf}\gf$ and let \rsmt\
and \mmt\ be the Ray-Singer and Milnor metrics on \dhb\
associated to the metrics \gft\ and \gtm.

By \cite[Theorem 6.3]{bz1}
\begin{multline}\label{E:def-lhs}
	\int_M\theta(F,\gft)(\nf)^*\mqf=
	     \int_M\thet(\nf)^*\mqf\\
			+2t\rk(F)\int_Mf\etm -2t\rk(F)\Tr_s^B[f].
\end{multline}
{}From  \refe{bz} and \refe{def-lhs}, we get
\begin{multline}\label{E:wit-milrs}
	\log\left(\frac\rsmt\mmt\right)^2=
		-\int_M\thet(\nf)^*\mqf  \\
			- 2t\rk(F)\int_Mf\etm+2t\rk(F)\Tr_s^B[f].
\end{multline}

\sec{mainres}{Proof of \refp{reduce} and \reft{mainth2}}

\ssec{1}{}
{}For each $t>0$, we equip \om\ with the inner product
\eq{in-pr-t}
	\langle\alp,\alp'\rangle_{\om,t}=
		\int_M\langle\alp\land\ast \alp'\rangle_{\gft},
\end{equation}
and we denote by  \rsmbt\  the
$L_2$ metric on \dhb\ (cf. \refss{badrsmet}) associated to this inner
product.

Let \mmbt\ be the Hodge-Milnor metric on \dhb\ (cf. \refd{mmb})
associated to the metric \gft\ on $F$.

{}From \refe{mm-mmb}, \refe{rsmet} and \refe{wit-milrs}, we get
\begin{multline}\label{E:pr1}
	\log \rstt - \log\mtt=\\
	  -\frac12\int_M\thet(\nf)^*\mqf-t\rk(F)\,\int_Mf\etm\\
		+t\rk(F) \Tr_s^B[f]+
			\log\left(\frac\mmbt\rsmbt\right).
\end{multline}
\refp{reduce} follows now from \refe{pr1} and the
following lemma:

\lem{m-rs} As $t\to+\infty$,
   \eq{lemma}
	\log \left(\frac{\mmbt}{\rsmbt}\right)=
        \left(\frac n4\chi(F)-\frac 12\chi'(F)\right)
        \log\left(\frac t\pi\right)+ o(1).
   \end{equation}
\elem
\prf
Let $d^{F\ast}_t$ $(t>0)$ be the formal adjoint of $d^F$
with respect to the inner product \refe{in-pr-t} and let
$\Del_t=d^Fd^{F\ast}_t+d^{F\ast}_td^F$.

Let \omot\  be the direct sum  of the eigenspaces
of $\Del_t$ associated to eigenvalues $\lam\in[0,1]$.
The pair $(\omot,d^F)$ is a subcomplex of $(\om,d^F)$ and
the inclusion induces an isomorphism of cohomology
\eq{h01=h}
	H^\b(\omot,d^F)\simeq H^\b(M,F).
\end{equation}

We denote by $\|\cdot\|_{\omot,t}$ the norm on \omot\ determined
by inner product \refe{in-pr-t}
and by $\|\cdot\|_{\ts,t}$ the norm on \ts\ determined
by \gft\ (cf. \refss{milmet}).

Recall that $\calF\in\End(\ts)$ was defined in \refss{wit-mt}.
Clearly,
\eq{pr2}
	\|\alp\|_{\ts,t}=\|e^{-t\calF}\alp\|_{\ts,0},
			\qquad (\alp\in\ts).
\end{equation}
{}For an operator $T:\ts\to \omot$ we denote by $T^*$ its adjoint
with respect to the norms $\|\cdot\|_{\ts,0}$ and
$\|\cdot\|_{\omot,t}$.

In the sequel, $\bold{o}(1)$ denotes an element of
$\End\big(\ts\big)$ which
preserves the $\ZZ$-grading and is $o(1)$ as $t\to\infty$.

By \cite[Theorem 6.9]{bz2}, if $t>0$ is large enough, there exists an
isomorphism
\[ e_t:\ts\to\omot \]
of $\ZZ$-graded  Euclidean vector spaces such that
\begin{equation}\label{E:et}
	e_t^*e_t=1+\bold{o}(1).
\end{equation}
By \cite[Theorem 6.11]{bz2}, for any $t>0$, there is
a quasi-isomorphism of complexes
\[ P_t:\Big(\omot,d^F\Big)\to \Big(\ts,\p\Big), \]
which induces the canonical isomorphism
\begin{equation}\label{E:isom}
	 \hb\simeq H^\b(\omot,d^F)\simeq H^\b(\ts,\p)
\end{equation}
and such that
\begin{equation}\label{E:pt-et}
    P_te_t=
     e^{t\calF}\left(\frac t\pi \right)^{n/4-N/2}\big(1+\bold{o}(1)\big).
\end{equation}
Here $e^{t\calF}\left(\frac{t}\pi\right)^{n/4-N/2}$
denotes the operator on \ts\ which, for $x\in B$, acts  on
$[W^u(x)]^*\otimes F_x$ by multiplication
by $e^{tf(x)}\left(\frac{t}\pi\right)^{n/4-index(x)/2}$.
In particular, for $t>0$
large enough, $P_t$ is  one to one.

Let $P_t^*$ denote the adjoint of $P_t$ with respect to the norms
$\|\cdot\|_{\ts,0}$ and $\|\cdot\|_{\omot,t}$.
{}From \refe{et}, \refe{pt-et} we get
\eq{pp*}
    P_tP_t^*=
     e^{2t\calF}\left(\frac{t}\pi\right)^{n/2-N}\big(1+\bold{o}(1)\big).
\end{equation}
Denote by $P_t^\#$ the adjoint of $P_t$ with respect to the norms
$\|\cdot\|_{\ts,t}$ and \text{$\|\cdot\|_{\omot,t}$}. Clearly,
\eq{P+}
	P_t^\#=P^*\cdot e^{-2t\calF}.
\end{equation}
Hence,
\eq{PP+}
	P_tP_t^\#=
	     \left(\frac t\pi\right)^{n/2-N}\big(1+\bold{o}(1)\big).
\end{equation}

{}Fix $0\le i\le n$. Let $\sig\in\h{i}$ and let $\ome_t\in \Ker \Del_t$
be the harmonic form representing $\sig$.

Denote by  $\p^\#_t$ the adjoint of \p\ with respect to the norms
$\|\cdot\|_{\ts,t}$ and \text{$\|\cdot\|_{\omot,t}$} and let
$\Pi:\ts\to\Ker(\p\p^\#_t+\p^\#_t\p)$ be the orthogonal projection.
Then $\Pi P_t\ome_t\in C^i(W^u,F)$ corresponds to $\sig$ via
the canonical isomorphisms \refe{hodge-ts}, \refe{isom}.

As $P_t$ commutes with $\p$, we see that
\eq{ker}
	P_t\ome_t\in \Ker\p,\qquad
	\left(\frac{t}\pi\right)^{n/2-i}\big(P_t^\#\big)^{-1}
	\ome_t\in\Ker \p^\#_t.
\end{equation}
By \refe{PP+}, we get
$\left(\frac{t}\pi\right)^{n/2-i}\big(P_t^\#\big)^{-1}\ome_t=
\big(1+\bold{o}(1)\big)P_t\ome_t$. Then \refe{ker} implies
\eq{pi1}
	\big\|\Pi P_t\ome_t\big\|_{\ts,t}=
	  \big\|P_t\ome_t\big\|_{\ts,t}\big(1+o(1)\big).
\end{equation}
It follows from \refe{PP+}, \refe{pi1} that
\eq{final}
	\big\|\Pi P_t\ome_t\big\|_{\ts,t}=
	  \left(\frac{t}\pi\right)^{n/4-i/2}
	    \big\|\ome_t\big\|_{\omot,t}\big(1+o(t)\big).
\end{equation}
By \refe{final} and by the definitions of the metrics
\mmbt, \rsmbt,  we obtain \refe{lemma}.
\epr
The proof of \refp{reduce} is completed.

\noindent
\reft{mainth2} follows now
from \refp{reduce}, \refl{mil-as} and \refp{mil-as-r}.


\end{document}